\documentclass[]{amsart}   

\usepackage{amsfonts}
\usepackage{amssymb}
\usepackage{mathrsfs}

\newtheorem{theorem}{Theorem}[section]
\newtheorem{lemma}[theorem]{Lemma}
\newtheorem{proposition}[theorem]{Proposition}
\newtheorem{corollary}[theorem]{Corollary}
\theoremstyle{definition}
\newtheorem{definition}[theorem]{Definition}
\newtheorem{example}[theorem]{Example}

\theoremstyle{remark}
\newtheorem{remark}[theorem]{Remark}
\numberwithin{equation}{section}

\newtheorem{problem}[theorem]{Problem}

\title[The uniform convergence topology on separable subsets]{The uniform convergence topology on separable subsets}

\date{June 30, 2023}

\author[J.\,A. Cruz-Chapital]{J.\,A. Cruz-Chapital$^*$}
\address{Centro de Ciencias Matemáticas, UNAM, Morelia}
\email{jorgeacruzchapital@ciencias.unam.mx}

\author[A.\,D. Rojas]{A.\,D. Rojas-S\'anchez} 
\address{Universidad Panamericana, Campus México, Cda. Augusto Rodin 498, Benito Juárez 03920, México}
\email{adrojas@up.edu.mx}

\author[\'A. Tamariz]{\'A. Tamariz-Mascar\'ua$^{**}$}
\address{Departamento de Matem\'aticas, Facultad de Ciencias, Circuito
  exterior s/n, Ciudad Universitaria, 04510, Ciudad de M\'exico, M\'exico}
\email{atamariz@unam.mx}

\author[H. Villegas]{H. Villegas-Rodr\'{\i}guez}
\address{Facultad de Ciencias F\'{\i}sico-Matem\'aticas, Universidad Aut\'onoma de Sinaloa, 
Ciudad Universitaria, Ave. de las Am\'ericas, Culiac\'an, Sinaloa, M\'exico}
\email{vrhi@uas.edu.mx}

\subjclass[2020]{54C35, 54A25, 54E52}  

\keywords{Function space, topology of uniform convergence, separable subsets,
  pseudouniform topology, Baire property, ordinal spaces, cellularity.}
  
\thanks{$*$ Part of the research of the first author was done 
while he participated in Thematic Program on Set Theoretic Methods in Algebra,
Dynamics and Geometry of the Fields Institute.}

\thanks{$**$ Corresponding author}
  
\begin{document}
\maketitle

\begin{abstract} 
For a topological space $X$, let $(\mathbb{R}^X)_s:= (\mathbb{R}^X,\mathcal{T}_s)$ be the cartesian product of $|X|$ copies of the real line $\mathbb{R}$ with the topology of the uniform convergence on separable subsets of $X$. In this article we analyze the subspace $C(X)$ of $(\mathbb{R}^X)_s$ of all real-valued continuous functions on $X$, denoted by $C_s(X)$. We determine when $C_s(X)$ is dense and when is closed in $(\mathbb{R}^X)_s$, and we obtain some results about the Baire property in $C_s(X)$. Finally, we determine the cellularity of $C_s([0,\alpha])$ where $[0,\alpha]$ is the space of ordinal numbers belonging to $\alpha + 1$ with its usual order topology. 
\end{abstract}

\section{Introduction, Notations and Basic Definitions}

The space $C_s(X)$ of all continuous real-valued functions provided with the topology of the uniform convergence on separable subsets of $X$ was first analized in \cite{PichTamVill} where several questions are posed. In this article, we continue with the study of this kind of topological spaces. First, we recall the constructions on $C(X)$ of uniform topologies on a semi-base of an ideal in $X$, and the $\alpha$-open topology where $\alpha$ is a $\pi$-network and their basic properties (Section 2). Moreover, we characterize the spaces $X$ for which $C_s(X)$ is dense in $(\mathbb{R}^X)_s$ (Section 3), analize the Baire property for spaces of the form $C_s(X)$ (Section 4), and in Section 5 we determine the cellularity of $C_s([0,\alpha))$ where $\alpha$ is an ordinal number.

\medskip 

All topological concepts that are not defined
here should be understood as in \cite{En}, and 
all topological spaces consider in this paper are assumed to be Tychonoff with more than one point. 

For spaces $X$ and $Y$ the symbol $X \cong Y$ means that $X$ is homeomorphic to $Y$.
When $X$ and $Y$ are two topological spaces with the same underlying
set, the symbol $X\leq Y$ means that the topology of $X$ is finer than
the topology of $Y$. When $X\leq Y$ and $X\not\cong Y$ we will write $X<Y$.  The closure of a subset $F$ of a topological space $X$ is denoted by ${\rm cl}_XF$, or briefly by ${\rm cl}F$. 

As usual, given a topological space $X$, $C(X)$ denotes the collection
of all continuous real-valued functions whose domain is $X$. If $f \in C(X)$, and $F \subseteq X$, $f \restriction F$ is the function belonging to $C(F)$ satisfying $(f \restriction F) (x) = f(x)$ for every $x \in F$.

Given a set $X$ and a cardinal $\kappa$, the symbol $[X]^{<\kappa}$
will denote the collection of all subsets of $X$ which have
cardinality $<\kappa$. A similar convention will apply to $[X]^{\leq
  \kappa}$. Finally, $[X]^\kappa$ is the set of all subsets of $X$
whose size is precisely $\kappa$.

The first infinite ordinal is $\omega$. The set $\omega\setminus \{0\}$ will be
denoted by $\mathbb N$. With $\omega_1$ we mean the first uncountable cardinal.
In what follows, $\lambda$, $\kappa$, $\tau$ and $\theta$ denote cardinal numbers, and $\alpha$, 
$\eta$, $\delta$, $\xi$, $\gamma$, $\mu$ and $\zeta$ denote ordinal numbers. For an ordinal number $\xi$, $[0,\xi)$ (resp., $[0,\xi]$) denotes the space of all ordinals strictly less than $\xi$ (resp., strictly less or equal to $\xi$) with its order topology, and for every set $X$, $D(X)$ is the discrete space of cardinality $|X|$.
The cardinality of the set of real numbers $\mathbb R$, will be
denoted by $\mathfrak c$. 

\section{Basic results on uniform convergence topologies}\label{sec:3}

In this section we recall the basic knowledge about uniform convergence topologies. The results presented here are known. The proofs can be found in \cite{MT}. 
In what follows, $\alpha$ is a non-empty collection of non-empty subsets of a topological space $X$.

Let $\mathcal{S}_\alpha$ be the collection of all sets of the form $[A,U] = \{g \in C(X) : g[A] \subseteq U\}$ where $A \in \alpha$ and $U$ is an open subset of 
$\mathbb{R}$. Let $\mathcal{B}_\alpha$ be the collection of all sets of the form $\bigcap \mathcal{S'}$ where $\mathcal{S'}$ is a finite subcollection of $\mathcal{S}$. Finally, let $\mathcal{T}_\alpha$ be the collection of all sets of the form $\bigcup \mathcal{B'}$ where $\mathcal{B'} \subseteq \mathcal{B}_\alpha$.

Now, let $\mathcal{B}_{\alpha,u}$ be the collection of all subsets of the form $V^\circ(f,A,\epsilon) = \{g \in C(X) : sup_{x \in A}|f(x) - g(x)| < \epsilon\}$ where $f \in C(X)$, $A \in \alpha$ and $\epsilon > 0$. And let 
$\mathcal{T}_{\alpha,u}$ be the collection of all subsets of the form $\bigcup \mathcal{B'}$ where $\mathcal{B}'$ is a subset of $\mathcal{B}_{\alpha,u}$.

\begin{proposition}\label{prop:1.1/2}
For each space $X$ and each $\alpha$, the collection $\mathcal{T}_\alpha$ is a topology on $C(X)$.
\end{proposition}

For each $\alpha$, $\mathcal{T}_\alpha$ is called the $\alpha$-open topology on $C(X)$. The set $C(X)$ equipped with this topology $\mathcal{T}_\alpha$ is denoted by $C_\alpha(X)$. 

\begin{definition}
A collection $\alpha$ of subsets of $X$ is an \textit{ideal} in $X$ if (i) $\emptyset \in \alpha$, (ii) $A \cup B \in \alpha$, whenever $A, B \in \alpha$, and (iii) every subset of an element in $\alpha$ belongs to $\alpha$. The collection $\alpha$ is a \textit{base of ideal} if for each $A,B \in \alpha$ there exists $C\in \alpha$ such that $A \cup B \subseteq C$. And 
$\alpha$ is a {\it semibase of an ideal} if for each $A,B \in \alpha$ there exists $C\in \alpha$ such that $A \cup B \subseteq \overline{C}$.
\end{definition}

\begin{proposition}\label{prop:1.1}
For each space $X$ and each $\alpha$, the following statements are equivalent:
\begin{enumerate}
\item the collection $\mathcal{T}_{\alpha,u}$ is a topology on $C(X)$;

\item the collection $\mathcal{T}_{\alpha,u}$ is a topology on $C(X)$ such that for each $f \in C(X)$, $\mathcal{V}^\circ(f) = \{V^\circ(f,A,\epsilon) : A \in \alpha, \epsilon > 0\}$ is a local base at $f$ of open subsets; 

\item $\alpha$ is a semibase of an ideal.
\end{enumerate}
\end{proposition}

For each semibase of an ideal $\alpha$, $\mathcal{T}_{\alpha,u}$ is the {\it uniform topology on the elements of $\alpha$ on $C(X)$}. 
Observe that in this topology, the collection of all sets of the form $V(f,A,\epsilon) = \{g \in C(X) : \forall x\in A (|f(x)-g(x)| < \epsilon)\}$ where $A \in \alpha$ and $\epsilon >0$ is a system of (not necessarily open) neighborhoods  of $C_{\alpha,u}(X)$ at $f$. We denote the topological space $(C(X),\mathcal{T}_{\alpha,u})$ by $C_{\alpha,u}(X)$. In the case where $X$ is an element in the semibase of ideal $\alpha$, then $\mathcal{T}_{\alpha,u}$ is denoted by $\mathcal{T}_u$ and is called the {\it uniform convergence topology on $X$}. In this case, $C_{\alpha,u}(X)$ is denoted by $C_u(X)$. It occurs that the collection of sets of the form $V^\circ(f,X,\epsilon)$ where $f \in C(X)$ and $\epsilon > 0$ form a base for $\mathcal{T}_u$.

\medskip

\noindent {\bf Notation:} Expressions such as $C_{\alpha,u}(X) \leq C_{\beta,u}(X)$, $C_\alpha(X) \leq C_\beta(X)$ and $C_{\alpha}(X) \leq C_{\beta,u}(X)$ mean 
$\mathcal{T}_{\alpha,u} \subseteq \mathcal{T}_{\beta,u}$, $\mathcal{T}_{\alpha} \subseteq \mathcal{T}_{\beta}$ and $\mathcal{T}_{\alpha} \subseteq \mathcal{T}_{\beta,u}$, respectively. 

\medskip

For each $\alpha$, we denote by $\widetilde{\alpha}$, $\overline{\alpha}$ and $I_{\widetilde{\alpha}}$ the collections $\{A \subseteq X : A$ is the union of a finite subcollection of $\alpha\}$, $\{\overline{A} : A \in \alpha\}$ and $\{B \subseteq X : \exists A \in \widetilde{\alpha} (B \subseteq A)\}$, respectively. 
Observe that for each $\alpha$, $\widetilde \alpha$ is a base of an ideal; moreover, $I_{\widetilde{\alpha}}$ is the ideal generated by $\widetilde{\alpha}$. 

\begin{proposition}\label{prop:1.2}
For each space $X$ and each $\alpha$, we have:
\begin{enumerate}
\item $\mathcal{T}_{\alpha} = \mathcal{T}_{\widetilde{\alpha}}$;
\item $\mathcal{T}_{\alpha} \subseteq \mathcal{T}_{I_{\widetilde{\alpha}}}$.
\end{enumerate}
\end{proposition}

A subset $A$ of a topological space $X$ is {\it $C$-compact} if for each $f \in C(X)$, $f[A]$ is a compact subset of $\mathbb{R}$.

\begin{proposition}\label{prop:equal}
Let $\alpha$ be a collection of $C$-compact subsets of $X$. Then, $\mathcal{T}_\alpha = \mathcal{T}_{\overline \alpha}$.
\end{proposition}

For each $\alpha$, we define $J_\alpha = \{B \subseteq X : B = cl_XB \text{ and } \exists A \in \alpha (B \subseteq A)\}$. 
Observe that if $\alpha$ is a semibase of an ideal, then $\overline \alpha$ and $J_{\overline {\alpha}}$ are bases of ideals. 

\begin{proposition}\label{prop:1.3}
For any space $X$ and any semibase of an ideal $\alpha$ of subsets of $X$ we have:
\begin{enumerate}
\item $\mathcal{T}_{\alpha,u} = \mathcal{T}_{\widetilde \alpha,u}$;
\item $\mathcal{T}_{\alpha,u} = \mathcal{T}_{\overline \alpha,u}$;
\item $\mathcal{T}_{\overline \alpha,u} = \mathcal{T}_{J_{\overline \alpha},u}$;
\item $\mathcal{T}_{\alpha,u} = \mathcal{T}_{I_{\overline \alpha},u} = \mathcal{T}_{\overline {\alpha},u} = \mathcal{T}_{I_{\alpha},u}$.
\end{enumerate}
\end{proposition}

\begin{remark}\label{rem:trivial}
\begin{enumerate}
\item The collection $\{X\}$ is a semibase of an ideal.

\item If $\alpha$ is a semibase of an ideal, then $\beta = \alpha \cup \{X\}$ is a semibase of an ideal too and $\mathcal{T}_{\alpha,u} \leq \mathcal{T}_{\beta,u} = \mathcal{T}_{\{X\}} = \mathcal{T}_u$. 

\item If each element in $\overline{\alpha}$ belongs to $\alpha$, then $\mathcal{T}_{\overline{\alpha}} \subseteq \mathcal{T}_\alpha$.
\end{enumerate}
\end{remark}

A collection $\alpha$ of subsets of a topological space $X$ is {\it closed under closed countable subsets} if when $A \in \alpha$ and $F\subseteq A$ is countable and closed in $A$, then $F \in \alpha$.

\begin{proposition}\label{prop:alpha,u}
Let $\alpha$ be a semibase of an ideal which is closed under closed countable subsets. Then, the following assertions are equivalent:
\begin{enumerate}
\item Each element of $\alpha$ is $C$-compact in $X$; 
\item $C_{\alpha}(X) = C_{\overline \alpha}(X) \leq C_{\overline \alpha,u}(X) = C_{\alpha,u}(X) \leq C_u(X)$;
\item $C_{\alpha}(X) \leq C_{\alpha,u}(X)$;
\item $C_\alpha(X) \leq C_u(X)$.
\end{enumerate}
\end{proposition}

A collection $\alpha$ of nonempty subsets of $X$ is a \textit{$\pi$-network} for $X$ if every nonempty open subset of $X$ contains at least an element of $\alpha$.

\begin{proposition}\label{prop:1.4}
For each $X$ and each $\alpha$, $\mathcal{T}_{\alpha}$ is Hausdorff if and only if $\alpha$ is a $\pi$-network for $X$. 
\end{proposition}

\begin{proposition}\label{prop:1.5}
For each $X$ and each semibase of an ideal $\alpha$ on $X$, $\mathcal{T}_{\alpha,u}$ is Hausdorff if and only if $\bigcup \alpha$ is a dense subset of $X$.
\end{proposition}

Observe that $\mathcal{T}_{\alpha,u}$ is a topology defined by an uniformity, so $C_{\alpha,u}(X)$ is completely regular. Since we are only considered Tychonoff spaces in this paper, in what follows every base for an ideal on $X$ will be consider, without explicit mention, such that $\bigcup \alpha$ is dense in $X$.

\begin{definition}
  Let $\kappa$ be an infinite cardinal. If $X$ is a topological space
  and $\alpha=[X]^{<\kappa}$, then we define
  $C_{\kappa,u}(X):=C_{\alpha,u}(X)$. In particular,
  \begin{enumerate}
  \item  $C_{\omega_1,u}(X)$  will be denoted by $C_s(X)$ 
  and its topology will be called {\sl the topology of
  pseudouniform convergence}.
  \item When $\kappa=\omega$, we obtain
  $C_p(X)$, the space of continuous functions on $X$ equipped with the
  topology of pointwise convergence.
  \end{enumerate}
\end{definition}

Let $X$ be a topological space. If $\alpha$ is the collection of all
compact subsets of $X$, $C_{\alpha,u}(X)$ will be denoted by
$C_k(X)$. It is a consequence of \cite[Theorem 1.2.3]{KunMcCoy}
that if one considers the cartesian product $\mathbb R^X$ endowed with
the compact-open topology (see \cite[Section 3.4]{En}), then the relative
topology of $C(X)$ coincides with the topology of $C_k(X)$.

When $\alpha=\{X\}$, we let $C_u(X):=C_{\alpha,u}(X)$. It is
straightforward to verify that $C_u(X)$ is the subspace $C(X)$ of
$\mathbb R^X$ endowed with the topology of uniform convergence (see
\cite[Section 2.6]{En}). In particular, $C_u(X)$ is metrizable.

It should be clear that $C_{\alpha,u}(X)\leq C_u(X)$ for any $\alpha$.

\begin{definition}
  Let $\alpha$ be a semibase for an ideal on a topological space $X$. We
  say that $\alpha$ is {\sl pseudouniform on $X$}  (or simply {\sl pseudouniform}
  when there is no risk of confusion) if for every sequence
  $\langle f_n:n\in\omega\rangle$ in $C(X)$ and all $f\in C(X)$ we
  have that $\langle f_n : n \in \omega \rangle$ converges to $f$ in $C_{\alpha,u}(X)$ ($f_n\stackrel \alpha \to f$) implies that  
  $\langle f_n : n \in \omega \rangle$ converges to $f$ in $C_{u}(X)$ ($f_n\stackrel u \to f$).
\end{definition}

Note that if $\alpha$ is pseudouniform on $X$ and $C_{\alpha,u}(X)\leq
C_{\beta,u}(X)$, then $\beta$ is pseudouniform as well; $\alpha=[X]^{<\omega_1}$ is a pseudouniform base of ideal, and 
all pseudouniform base of an ideal have dense union, and then it defines a Tychonoff topology on $C(X)$ (see \cite{PichTamVill}). 

\begin{theorem}\label{ptv_pseudo-cellular}
  If $\alpha$ is a semibase for an ideal on $X$ with dense union, the following
  statements are equivalent.
  \begin{enumerate}
  \item $\alpha$ is pseudouniform.
  \item For every
  family $\{U_n:n\in \omega\}$ of nonempty open sets in $X$ there exists $A\in \alpha$ such
  that the set $\{n\in\omega:U_n\cap A\neq\emptyset\}$ is infinite.
\item For every cellular 
  family $\{U_n:n\in \omega\}$ of non empty open sets of $X$ there exists $A\in \alpha$ such
  that the set $\{n\in\omega:U_n\cap A\neq\emptyset\}$ is infinite.
  \end{enumerate}
\end{theorem}

\begin{example}(Example of a topological space $X$ such that $C_k(X)$ is a pseudouniform topology satisfying $C_k(X) < C_s(X) < C_u(X)$.)
Let $\mathcal{A}$ be a maximal collection of subsets of $\omega_1$ which satisfies the following conditions: (i) Each element in $\mathcal{A}$ is countable and infinite, (ii) if $A, B \in \mathcal{A}$ are different, then $|A \cap B| < \omega$. Let $X= \psi(A)$ be the Mr\'owka-Isbell space on $\omega_1$. Then, $X$ is pseudocompact and locally compact. 
Thus $X$ is almost pseudo-$\omega$-bounded (Proposition 4.9 en \cite{PichTamVill}). This means that $C_k(X)$ is pseudouniform. Moreover, every compact subset of $X$ is separable. This fact implies that $C_k(X) \leq C_s(X)$ (Proposition 5.3.(2) in \cite{PichTamVill}). Since $X$ is not separable, $C_s(X) \not= C_u(X)$. 

Now we are going to verify the inequality $C_k(X) < C_s(X)$. If $C_k(X) = C_s(X)$ then $X$ would be $\omega$-bounded (Proposition 5.3.(1) in \cite{PichTamVill}). But each infinite countable subcollection of $\mathcal{A}$ is closed and discrete, so $X$ is not $\omega$-bounded ({\sl $\omega$-bounded} = the closure of each countable subset is compact). Therefore, $C_k(X)$ is pseudouniform and $C_k(X) < C_s(X) < C_u(X)$. In particular, $C_s(X)$ is not minimal among the pseudouniforme topologies of $C(X)$.
\end{example}

\begin{example}(Example of a topological space $X$ such that $C_k(X)$ is a pseudouniform topology with $C_s(X) < C_k(X) < C_u(X)$.)
Let $\kappa$ be a cardinal number $> \omega_1$ and of uncountable cofinality. Let $X = [0,\kappa)$. Since $X$ is pseudocompact and locally compact, then $X$ is almost pseudo-$\omega$-bounded. (Even more, $[0,\kappa)$ is almost pseudo-$\omega$-bounded iff $[0,\kappa)$ is countably compact iff ${\rm cof}(\kappa) > \omega$.) That is $C_k(X)$ es pseudouniform. Furthermore, $X$ is 
$\omega$-bounded  because ${\rm cof}(\kappa) > \omega$. Thus, $C_s(X) \leq C_k(X)$. On the other hand, $C_s(X) < C_k(X)$ because if we assume the contrary, every compact subset of $X$ would be contained in a separable subset of $X$ (Proposition 5.3 in \cite{PichTamVill}), but $[0,\omega_1]$ is compact and it is not contained in a separable subspace of $[0,\kappa)$ (The closure of every countable subset of $[0,\kappa)$ is still countable). Then, $C_s(X) < C_k(X)$. Moreover, $C_k(X) < C_u(X)$ because $X$ is not  compact \cite[Theorem 1.1.4, p. 5]{McNt}. (By the way, $C_k([0,\kappa))$ is a $k$-espace iff ${\rm cof}(\kappa) = \omega$. See Theorem 4.7.4, p. 64, and exercise 11,  page 71 in \cite{McNt}.)
\end{example}

\section{When is $C_s(X)$ dense in $(\mathbb{R}^X)_s$?}

\begin{definition}
A space $X$ is $C_\omega$-discrete if each countable subset of $X$ is $C$-embedded in $X$.
\end{definition}

Observe that if $X$ is $C_\omega$-discrete, then every countable subset of $X$ is closed and discrete. 

\begin{lemma}\label{prop:denso}
If $C_s(X)$ is dense in $(\mathbb{R}^X)_s$, then for each countable subset $N = \{x_n : n < \omega\}$ of $X$ there exists a collection $\mathcal{V} = \{V_n : n < \omega\}$ such that $V_n$ is a closed neighborhood of $x_n$ in $X$ for each $n < \omega$, such that the elements of $\mathcal{V}$ are pairwise disjoint and $\mathcal{V}$ is locally finite in $X$. In particular, $N$ is closed and discrete.
\end{lemma}

\begin{proof}
If $N$ is finite, the conclution is obvious. Assume then that $x_n \not= x_m$ if $n \not= m$. Consider the function $f \in \mathbb{R}^X$ defined as $f(x_n) = n$ and $f(y) = 0$ if $y \not\in N$. Since $C_s(X)$ is dense in $(\mathbb{R}^X)_s$, then there exists $g \in C(X)$ such that $g \in \widetilde{V}(f;N;1/4) = \{h \in \mathbb{R}^X : \forall x \in N, |f(x) - h(x)| < 1/4\}$. It occurs that $x_n \in B_n := g^{-1}[(n - 1/4,n + 1/4)]$. Let $V_n = g^{-1}[[n - 1/4,n + 1/4]]$. Each $V_n$ is a closed neighborhood of $x_n$ in $X$. We prove now that $\mathcal{V} = \{V_n : n < \omega\}$ is a locally finite family in $X$. Let $a \in X$ and let $g(a) = r \in \mathbb{R}$. The collection $\mathcal{A} = \{[n - 1/4,n+1/4] : n < \omega\}$ is locally finite in $\mathbb{R}$. Therefore, there is  a neighborhood $C$ of $r$ which intersect only a finite subcollection of elements in $\mathcal{A}$. The set $g^{-1}[C]$ is a neighborhood of $a$ which intersect only a finite subcollection of elements in $\mathcal{V}$. 
\end{proof}

The following result shows that the concept of pseudouniformly dense introduced in \cite{RTV} coincide with the $C_\omega$-discreteness. 

\begin{theorem}\label{theo:denso}
$C_s(X)$ is dense in $(\mathbb{R}^X)_s$ if and only if $X$ is $C_\omega$-discrete. 
\end{theorem}

\begin{proof}
Let $X$ be $C_\omega$-discrete. Let $f \in \mathbb{R}^X$ and let $\widetilde{V}(f;N;\epsilon) = \{h \in \mathbb{R}^X : \forall x \in N, |f(x) - h(x)| < \epsilon\}$ be a canonical neighborhood of $f$ in $(\mathbb{R}^X)_s$ where $N = \{x_n : n < \omega\}$ with $x_n \not= x_m$ if $n \not= m$. The set $N$ is discrete, so $f\restriction N : N \to \mathbb{R}$ is continuous. Since $X$ is $C_\omega$-discrete, then $N$ is $C$-embedded in $X$. This means that there is a function $h \in C(X)$ such that $h\restriction N = f$. Hence 
$h \in \widetilde{V}(f;N;\epsilon)$.

\medskip

Now assume that $C_s(X)$ is dense in $(\mathbb{R}^X)_s$. Let $N \subseteq \mathbb{R}$ be equal to $\{x_n : n < \omega\}$ with $x_n \not= x_m$ if $n \not= m$. By Lemma \ref{prop:denso} we find, for each $n < \omega$, an open neighborhood $B_n$ of $x_n$ in $X$ such that $\{cl_XB_n : n < \omega\}$ is a locally finite collection in $X$ which elements are pairwise disjoint. Let $g : N \to \mathbb{R}$ be a function; since $N$ is discrete (Lemma \ref{prop:denso}), $g$ is continuous. For each $n < \omega$, let $h_n : X \to \mathbb{R}$ be a continuous function such that $h_n(x_n) = g(x_n)$ and $h_n(y) = 0$ if $y \not\in B_n$ ($X$ is Tychonoff). Since $\{B_n : n < \omega\}$ is a pairwise disjoint locally finite collection, the function $h = \Sigma_{n<\omega} h_n$ is a well defined continuous function. Moreover, $h \restriction N = g$. That is, $N$ is $C$-embedded in $X$.   
\end{proof}

\begin{corollary}
If $X$ is pseudocompact and infinite, then $C_s(X)$ is not dense in $(\mathbb{R}^X)_s$.
\end{corollary}

\begin{proof}
Indeed, if $X$ is pseudocompact, it does not have a $C$-embedded copy of $\mathbb{N}$ (\cite[Corollary 1.20]{GJ}), so $X$ is not $C_\omega$-discrete.
\end{proof}

\begin{example}
There are spaces $X$ in which the closure of every countable subset is $C$-embedded in $X$, but $C_s(X)$ is not dense in $(R^X)_s$. This is the case of every infinite compact space $X$. 
\end{example}

\begin{lemma}\label{lema:uniforme}
Let $u$ be an Euclidean and bounded uniformity on $\mathbb{R}$. If $N$ is a countable topological space, then 
$(\mathbb{R}^N)_s = (\mathbb{R}^N)_u$.
\end{lemma}

\begin{proof}
Let $f \in \mathbb{R}^N$ and let $V^\circ(f;N',\epsilon)= \{\xi \in \mathbb{R}^N : sup_{x\in N'}|f(x) - \xi(x)| < \epsilon\}$ a canonical neighborhood of $f$ in $(\mathbb{R}^N)_s$. Let $g \in V(f;N',\epsilon)$. It follows that there exists $\delta < \epsilon$ such that $|f(x) - g(x)| \leq \delta$ for each $x \in N'$. Then $g \in V(g;N,(\epsilon - \delta)/2) \subseteq V(f;N',(\epsilon -\delta)/2n)$ because $h \in V(g;N,(\epsilon - \delta)/2)$, then for each $x \in N'$, $|f(x) - h(x)| \leq |f(x) - g(x)| + |g(x)-h(x)| \leq \delta + (\epsilon - \delta)/2 = (\delta + \epsilon)/2 < \epsilon$. This prove that $\mathcal{T}_s \subseteq \mathcal{T}_{u}$ in $\mathbb{R}^N$. 

Now, for each $f \in \mathbb{R}^N$, $V^\circ (f;N;\epsilon)$ is an element of $\mathcal{T}_s$ in $\mathbb{R}^N$ because $N$ is countable.
\end{proof}

\begin{theorem}\label{theo:cl}
Let $X$ be a topological space. If $f \in \mathbb{R}^X$ belongs to $cl_{(\mathbb{R}^X)_s}C_s(X)$, then 
for every $N \in [X]^{\leq \omega}$, $f\restriction N$ is continuous. 
\end{theorem}

\begin{proof}
Let $f \in cl_{(\mathbb{R}^X)_s}C_s(X)$ and let $N \in [X]^{\leq \omega}$. A local base of neighborhoods of $f$ in $(\mathbb{R}^X)_s$ is constituted by the collection $\{\widetilde{V}(f;N;\epsilon) : N \in [X]^{\leq \omega}\}$ where $\widetilde{V}(f;N;\epsilon) = \{g \in \mathbb{R}^X : \forall x \in N, |f(x)-g(x)| < \epsilon\}$. So, since $f \in cl_{(\mathbb{R}^X)_s}C_s(X)$, for each $\epsilon > 0$, $\widetilde{V}(f;N;\epsilon)$ contains a function $g \in C_s(X)$. In particular, each uniform neighborhood $\widetilde{V}(f\restriction N;\epsilon)$ of $f \restriction N$ in $(\mathbb{R}^N)_s$ contains a continuous function (in this case $g\restriction N$). That is, $f\restriction N$ belongs to $cl_{(\mathbb{R}^N)_s}C_s(N)$. But $s$ is the topology 
of the uniform convergence on $\mathbb{R}^{N}$ because $N$ is countable; that is $C_s(N) = C_u(N)$, where $u$ is the uniformity on $\mathbb{R}$ generated by a bounded metric 
equivalent to the Euclidean uniformity. But $(\mathbb{R}^{N})_u$ is metrizable (see Theorem 4.2.20 in \cite{En}), and so, there exists a sequence $(f_n)_{n < \omega}$ of elements 
in $C_u(N)$ which uniformely converges to $f\restriction N$. By Theorem 4.2.19, p. 264, in \cite{En}, this means that $f\restriction N$ is continuous. 
\end{proof}

\begin{proposition}\label{prop:cl1}
Let $X$ be a topological space and let $f \in \mathbb{R}^X$. If $f$ belongs to $cl_{(\mathbb{R}^X)_s}C_s(X)$, then 
for each $N \in [X]^{\leq \omega}$, $f\restriction N \in cl_{C_s(N)}\pi_N[C_s(X)]$. 
\end{proposition}

\begin{proof}
Assume that $f \in \mathbb{R}^X$ belongs to $cl_{(\mathbb{R}^X)_s}C_s(X)$. Let us see that if $N \in [X]^{\leq \omega}$, then $f\restriction N \in cl_{C_s(N)}\pi_N[C_s(X)]$. Let $N' \in [N]^{\leq \omega}$ and let $\epsilon >0$. Thus, there exists $g \in C_s(X)$ such that $g \in \widetilde V(f;N';\epsilon)$. It follows that $g \restriction N \in C_s(N)$ and $\widetilde V(f;N';\epsilon) \cap C_s(N) = \{h \in C(N) : \forall x \in N', |h(x) - f\restriction N(x)| < \epsilon\} = V(f\restriction N; N';\epsilon)$ whih is a canonical neighborhood of $f$ in $C_s(N)$. This means that $f\restriction N \in cl_{C_s(N)}\pi_N[C_s(X)]$.
\end{proof}

For each topological space $X$ and each $x$ in $X$, we denote by $\mathcal{V}(x)$ the collection of all neighborhoods of $x$. 

\begin{proposition}\label{prop:tight}
Let $X$ be a space with countable tightness ($t(X) = \omega$). If $f \in \mathbb{R}^X$ and for every $N \in [X]^{\leq \omega}$, $f\restriction N$ is continuous in $N$, then $f$ is continuous  in $X$.
\end{proposition}

\begin{proof}
Assume that $f$ satisfices the conditions of the proposition and $f$ is discontinuous at $z \in X$. Then, there exists $\epsilon > 0$ such that for every $V \in \mathcal{V}(z)$ there is $x_V \in V$ such that $f(x_V) \not\in (f(z)-\epsilon,f(z)+\epsilon)$. Of course, in particular $x_V \not= z$. Then $z \in cl_X\{x_V : V \in \mathcal{V}(z)\}$. Since $t(X) = \omega$, there exists $N \subseteq \{x_V : V \in \mathcal{V}(z)\}$, countable, such that $z \in cl_XN \setminus N$. By hypothesis, it occurs that $f\restriction N \cup \{z\}$ is continuous. On the other hand there exists a neighborhood $B$ of $z$ in $N\cup \{z\}$ such $f(B) \subseteq (f(z)-\epsilon,f(z) + \epsilon)$. This means that several $y \in N$ satysfy $f(y) \in (f(z)-\epsilon,f(x)+\epsilon)$. But $N \subseteq \{x_V : V \in \mathcal{V}(z)\}$, and we obtain a contradiction. So, $f$ must be continuous in $X$. 
\end{proof}

\begin{proposition}\label{prop:tight1}
If $X$ has countable tightnes, then $C_s(X)$ is closed in $(\mathbb{R}^X)_s$.
\end{proposition}

\begin{proof}
Indeed, if $f \in cl_{(\mathbb{R}^X)_s}C_s(X)$, then for every $N \in [X]^{\leq \omega}$, $f\restriction N \in C(N)$ (Theorem \ref{theo:cl}). By Proposition \ref{prop:tight}, $f \in C(X)$. 
\end{proof}

\begin{proposition}
The following statements are equivalent for a space $X$:
\begin{enumerate}
\item $C_s(X)$ is closed in $(\mathbb{R}^X)_s$; 
\item $f \in C(X)$ if and only if for each $N \in [X]^{\leq \omega}$ and each $\epsilon > 0$, there exists $g \in C(X)$ such that for any $x \in N$, $|f(x) - g(x)| < \epsilon$;
\item $f \in C(X)$ if and only if for each $N \in [X]^{\leq \omega}$ and each $\epsilon > 0$, there exists $g \in C(X)$ such that $g \in V(f,N,\epsilon)$. 
\end{enumerate} 
\end{proposition}

\begin{theorem}\label{theo:numerable}
Let $X$ be a space with countable tightness and $f \in \mathbb{R}^X$. Then, the following conditions are equivalent:
\begin{enumerate} 
\item $f$ is continuous; 
\item for every $N \in [X]^{\leq \omega}$, $f\restriction N$ is continuous; 
\item for every $N \in [X]^{\leq \omega}$, $f \restriction cl_XN$ is continuous. 
\end{enumerate}
\end{theorem}

\begin{proof}
The implications (1) $\Rightarrow$ (3) and (3) $\Rightarrow$ (2) are obvious. The implication (2) $\Rightarrow$ (1) is Proposition \ref{prop:tight}.
\end{proof}

\begin{proposition}\label{lem:discreto}
$C_s(X) = (\mathbb{R}^X)_s$ if and only if $X$ is discrete.
\end{proposition}

\begin{proof}
If $X$ is discrete, then $C(X) = \mathbb{R}^X$. Thus, $C_s(X) = (\mathbb{R}^X)_s$. Now, if $C_s(X) = (\mathbb{R}^X)_s$, then every function defined on $X$ is continuous. So, if $x \in X$ the characteristic function $\chi_{\{x\}}$ is continuous. This implies that $\{x\}$ is open. So $X$ is discrete. 
\end{proof}

\begin{corollary}
A space $X$ is discrete if and only if $X$ is $C_\omega$-discrete and has countable tightness.
\end{corollary}

\begin{proof}
Of course, if $X$ is discrete, then it is $C_\omega$-discrete and with countable tightness. Now, since $X$ is $C_\omega$-discrete, then $C_s(X)$ is dense in $(\mathbb{R}^X)_s$ (Theorem \ref{theo:denso}) and since $X$ has countable tightness, $C_s(X)$ is closed in $(\mathbb{R}^X)_s$ (Proposition \ref{prop:tight1}). So, $C_s(X) = (\mathbb{R}^X)_s$. By Proposition  \ref{lem:discreto}, $X$ is discrete.
\end{proof}

\begin{proposition}\label{prop:tightness}
If $X$ has countable tightness, then $A \subseteq X$ is closed in $X$ if and only if $A \cap N$ is closed in $N$ for every $N \in [X]^{\leq 
\omega}$.
\end{proposition}

\begin{proof}
Suppose that $X$ have countable tightness. Let $A \subseteq X$. Obviously, if $A$ is closed in $X$, then for any $N \in [X]^{\leq 
\omega}$, $A \cap N$ is closed in $N$. Suppose then that $A \cap N $ is closed in $N$ for all $N \in [X]^{\leq \omega}$. We are going to show that $A$ is closed in $X$. Let $x \in cl_XA$. Since $X$ has countable tightness, there exists $B \subseteq A$ countable such that $x \in cl_XB$. Then $B \cup \{x\}$ is a countable set. Therefore, 
$(B \cup \{x \}) \cap A$ is closed in $B \cup \{x\}$. If $x \not \in A$, then $ (B \cup \{x\}) \cap A = B$ is closed in $B \cup \{x\}$. This means that there is an open $U$ in $X$ such that $x \in U$ and $B \cap U = \emptyset$, which is a contradiction to the hypothesis 
$x \in cl_XB$. So $x \in A$. This means that $A$ is a closed subset of $X$.
\end{proof}

\section{Completeness type properties in $C_s(X)$}

A space is called \emph{\v{C}ech-complete} if it is a set of type 
$G_{\delta }$ in some (hence in any) Hausdorff compactification of it. For
metrizable spaces, \v{C}ech completeness is equivalent to metrizability by a
complete metric, so $C_{u}(X)$ is \v{C}ech-complete, since it is a complete
metric space. It is well known that every \v{C}ech-complete space is a $k$-space.
So, by Theorem 5.4 in \cite{RTV} we have:

\begin{proposition}
$C_{s}(X)$ is a \v{C}ech-complete space if and only if $d(X)=\omega$.
\end{proposition}

A space $X$ is a \emph{Baire space} if for every sequence $G_{1},G_{2},...$
of open dense subset of $X$, the intersection $\bigcap_{i=1}^{\infty }G_{i}$ is a dense set. Every \v{C}ech-complete space is a Baire
space \cite{En}.

\begin{theorem}
If $C_s(X)$ is a Baire space, then every pseudocompact subspace of $X$ is contained in a separable closed subset of $X$.
\end{theorem}

\begin{proof}
Assume that $A \subseteq X$ is pseudocompact, and that for every separable closed subset $F$ of $X$, $A \setminus F \not= \emptyset$.
Put, for each $n < \omega$, $G_{n}=\{f\in C_{s}(X): \exists x \in A(f(x) > n)\}$.
We show that $G_{n}$ is open and dense in $C_{s}(X)$. Let $f\in G_{n}$. Take an $x\in
A$ such that $f(x)>n$. Put $\delta =f(x)-n>0$ and consider the canonical 
neighborhood of $f$,  $V(f,x,\delta )$ in $C_{s}(X)$. Clearly, $f\in
V(f,x,\delta )\subseteq G_{n}$. Hence $G_{n}$ is open in $C_{s}(X)$. 

Now, let $V(g,E,\epsilon )$ be a canonical neighborhood of $g$ in $C_{s}(X)$. 
By hypothesis, there is a $y\in A\setminus \overline{E}$. 
Since $X$ is Tychonoff there is  
$f\in C_{s}(X)$ such that $f(y) = n + 1 - g(y)$ and $f(z) = 0$ for every $z \in \overline{E}$.
Define $h = f + g$. Then $h\in G_{n}\cap V(g,E,\epsilon )$. 
Hence $\overline{G_{n}}=C_{s}(X)$. 

Since $C_s(X)$ is a Baire space, there is  $l \in \cap
\{G_{n}:n<\omega \}$. Then for each $n<\omega$ there is a $x_{n}\in A$ 
such that $l(x_{n})>n$, which is impossible since $A$ is pseudocompact.
\end{proof}

\begin{corollary}
Let $X$ be a non separable space. If $C_s(X)$ is a Baire space, then $X$ is not pseudocompact. 
\end{corollary}

\begin{theorem} 
If $X$ is $C_\omega$-discrete, then $C_s(X)$ (resp., $(\mathbb{R}^X)_s$) is a Baire space.
\end{theorem}

\begin{proof}
Let $\langle D_{n}:n<\omega \rangle$ be a sequence of open dense subsets of $C_{s}(X)$. Let $U$ be an arbitrary non-empty open set in $C_{s}(X)$. 
We define four sequences by induction: a sequence 
$\langle f_{n}:n<\omega \rangle$ of elements in $C_{s}(X)$; a sequence $\langle E_{n}:n<\omega\rangle$ of countable subsets of $X$; a sequence $\langle W_{n}:n<\omega
\rangle$ of open subsets in $C_{s}(X)$; and a sequence $\langle \epsilon _{n}:n<\omega\rangle$
of real numbers; all these satisfying, for each $n<\omega$, the conditions:

\begin{enumerate}
\item[(a)] $V^\circ(f_{n},E_{n},\epsilon _{n})\subseteq W_{n}\subseteq cl_{C_s(X)}{W_{n}} 
\subseteq D_{n} \cap U$,

\item[(b)] $V^\circ(f_{n+1},E_{n+1},\epsilon _{n+1})\subseteq V^\circ(f_{n},E_{n},\epsilon
_{n})$,

\item[(c)] $E_{n}\subseteq E_{n+1}$, 

\item[(d)] $\epsilon _{n+1} < \epsilon _{n} < 1/2^{n+1}.$
\end{enumerate}

We choose an arbitrary function $f_{0}\in D_{0} \cap U.\ $Since $C_{s}(X)$ is a
regular space, there is an open set $W_{0}$ of $C_s(X)$ such that\ $f_{0}\in W_{0}\subseteq 
cl_{C_s(X)}(W_{0})\subseteq D_{0}\cap U$; and by definition of the canonical neighborhoods of 
$C_{s}(X)$, there is a countable set $E_{0}\subseteq X$, and $\epsilon
_{0}>0$, that we can suppose to be $<\frac{1}{2}$, such that $V^\circ(f_{0},E_{0},\epsilon
_{0})\subseteq V_{0}$. Assume that we have already choose $f_0,\dots,f_n$, $W_0,\dots,W_n$, $E_0,\dots,E_n$ and $\epsilon_0,\dots,\epsilon_n$ satisfying conditions (a)-(d).   

Since $D_{n+1}$ is dense, $V^\circ(f_{n},E_{n},\epsilon
_{n})\cap D_{n+1} \cap U \neq \varnothing$. We choose then a function 
$f_{n+1}\in V^\circ (f_{n},E_{n},\epsilon _{n})\cap D_{n+1} \cap U$. Repeating the
argument for the step 0, we can chose an open set $W_{n+1}$, a countable set $E_{n+1}\subseteq X$, 
and $\epsilon _{n+1}>0$ such that $V^\circ(f_{n+1},E_{n+1},\epsilon
_{n+1})\subseteq W_{n+1}\subseteq cl_{C_s(X)}W_{n+1}\subseteq V^\circ(f_{n},E_{n},\epsilon
_{n})\cap D_{n+1} \cap U$. Notice that $E_{n+1}$ and $\epsilon_{n+1}$ can be choosen such that 
$E_n \subseteq E_{n+1}$ and $\epsilon _{n+1}<min\{\epsilon _{n},1/2^{n+2}\}$. Continuing
this process, we obtain the required sequences.

\noindent {\bf ($\star$)} \hskip .25cm Observe that for every $n > k$ and for every $x \in E_k$, $|f_n(x) - f_k(x)| < \epsilon _k$, because $f_n \in V^\circ(f_n,E_n,\epsilon_n) \subseteq V^\circ(f_k,E_k,\epsilon_k)$. 

Let $E$ be equal to $\bigcup_{n<\omega}E_n$. Since $X$ is $C_\omega$-discrete, $E$ is discrete. We define now $\langle g_{n}:n<\omega \rangle$ in $C(E)$ by 
\[
g_{n}(x)=\left \{ 
\begin{tabular}{l}
$f_{n}(x)$\ if\ $x\in E_{n}$ \\ 
$f_{n+m}(x)$ if $x\in E_{n+m}\backslash E_{n+m-1},\ m=1,2,3,....$ 
\end{tabular}%
\right. 
\]%

The sequence $\langle g_{n}:n<\omega \rangle$ is Cauchy in $C_u(E)$. 
Indeed, let $\epsilon >0$ and let $N_\epsilon < \omega$ such that $\Sigma_{i\geq N_\epsilon}(1/2^{i}) < \epsilon$. For each $i \in \mathbb{N}$, if $x \in E_i$, $|g_i(x) - g_{i+1}(x)| = |f_i(x) - f_{i+1}(x)| < 1/2^{i+1}$, and if $i \in E \setminus E_i$, $g_i(x) = g_{i+1}(x)$. Then, if $k > m \geq N_\epsilon$, $|g_m(x) - g_k(x)| < \Sigma_{i=m+1}^{i=k} (1/2^{i}) \leq \Sigma_{i\geq N_\epsilon}(1/2^{i}) < \epsilon$. 

By Lemma \ref{lema:uniforme} and because $C_u(E)$ is completely metrizable (indeed, the topology of $C_u(X)$ is defined by the bounded complete metric $\rho$ defined on $\mathbb{R}$ by $\rho(x,y) = min\{|x-y|, 1\}$, and the uniformity defined by $\rho$ in $C_(X)$ and that defined by the absolute value are the same; see Theorem 4.2.20 in \cite{En} and Theorem 1.2.6 in \cite{McNt}), $\langle g_{n}:n<\omega \rangle$ converges uniformly to some (continuous) function $g\in C(E)$. 

Since $X$ is $C_\omega$-discrete and $E$ is countable, there is a continuous function $h \in C(X)$ such that $h\restriction E = g$. Finally, we are going to show that $h\in D_{k}\cap U$ for each $k<\omega$. Let $k < \omega$ be fixed. 

Since $\langle g_n : n < \omega \rangle$ uniformely converge to $g$ in $E$, for each $\epsilon > 0$, there is $N_\epsilon < \omega$ such that $|g(x) - g_n(x)| < \epsilon$ for every $x \in E$ and every $n \geq N_\epsilon$.
So, for every $x \in E_k$ and every $n \geq max\{N_\epsilon,k\}$, 
$|g(x) - f_k(x)| \leq |g(x) - f_n(x)| + |f_n(x) - f_k(x)| = |g(x) - g_n(x)| + |f_n(x) - f_k(x)| < \epsilon + |f_n(x) - f_k(x)|$. Because ($\star$), above, we conclude that $|g(x) - f_k(x)| < \epsilon + \epsilon _k$ for every $x \in E_k$. 
Since $\epsilon$ is an arbitrary positive real number, we have therefore the following relations  $|g(x) - f_k(x)| = |h(x) - f_k(x)| \leq \epsilon_k$ for every $x \in E_k$. But this implies that $h \in cl_{C_s(X)}V^\circ(f_k,E_k,\epsilon_k)$.  

Indeed, let $V(h,H,\epsilon)$ be a canonical neighborhood of $h$. First case: $H \cap E_k = \emptyset$. We have that $H \cup E_k$ is countable and so discrete. Then, the relation $t : H \cup E_k \to \mathbb{R}$ defined by $t(x) = h(x)$ if $x \in H$ and $t(x) = f_k(x)$ if $x \in E_k$ is a function (and then a continuous one on $H \cup E_k$). Since $X$ is $C_\omega$-discrete, there is a continuous map $t' \in C(X)$ such that $t' \restriction H \cup E_k = t$. The function $t'$ belongs to $V(h,H,\epsilon) \cap V(f_k,E_k,\epsilon_k)$. Second case: $H \cap E_k = M$ is not empty. Again, $H \cup E_k$ is discrete being countable. Let $\epsilon > 0$. Take $\delta > 0$ such that $\delta < min\{\epsilon,\epsilon_k\}$. We define the relation $t : H \cup E_k \to \mathbb{R}$ as follows: 
$t(x) = h(x)$ if either $x \in H \setminus M$ or $x \in M$ and $|f_k(x) - h(x)| < \epsilon_k$, and $t(x)= f_k(x)$ if $x \in E_k \setminus M$. If $x \in M$ and $|f_k(x) - h(x)| = \epsilon_k$, we have two possible cases: First, 
$h(x) = f_k(x) - \epsilon_k$; in this case we define $t(x) = h(x) + \delta$. Second case: $h(x) = f_k(x) + \epsilon_k$; in this case, we define $t(x) = h(x) - \delta$. The relation $t$ belongs to $C(H\cup E_k)$. Since $X$ is $C_\omega$-discrete, there is $t' \in C(X)$ such that $t' \restriction (H \cup E_k) = t$. In any case we have that $t' \in V(h,H,\epsilon) \cap V(f_k,E_k,\epsilon_k)$. All this proves that $h \in cl_{C_s(X)}V(f_k,E_k,\epsilon_k)$.  

\medskip

Since $cl_{C_s(X)}V^\circ(f_k,E_k,\epsilon_k) \subseteq cl_{C_s(X)}W_k \subseteq D_k \cap U$, we have that
$h \in (\bigcap_{n<\omega}D_k)\cap U$; that is, $\bigcap_{n<\omega}D_k$ is a dense subset of $C_{s}(X)$. 
Therefore, $C_s(X)$ is a Baire space.

On the other hand, if $C_s(X)$ is a Baire space, then $(\mathbb{R}^X)_s$ is too because $C_s(X)$ is dense in $(\mathbb{R}^X)_s$. 
\end{proof}

The following concepts appears in \cite{McNt}. A subfamily $\mathcal{F}$ of a family $\mathcal{G}$ will be said to {\it move off} $\mathcal{G}$ provided that for each $G \in \mathcal{G}$, there exists $F \in \mathcal{F}$ which is disjoint from $G$. A family of subsets of $X$ is {\it strongly discrete} provided that the members of the family have neighborhoods which form a discrete family in $X$.

\begin{theorem}
Let $X$ be a non-separable $\omega$-normal topological space. If $C_s(X)$ is a Baire space, then each move off family of separable closed subsets of $X$ contains a strongly discrete countable collection. 
\end{theorem}

\begin{proof}
For each $n < \omega$, let $D_n = \bigcup_{E\in \mathcal{F}}{\rm int}_{C_s(X)}V(n+1/2,\overline E,1/2)$, where $\mathcal{F}$ moves off the family of separable closed subsets of $X$. Of course $D_n$ is open. We claim that $D_n$ is dense in $C_s(X)$. Indeed, let $f \in C_s(X)$ and let $V(f,H,\epsilon)$ be a canonical neighborhood of $f$. Since $\mathcal{F}$ moves off the collection of all separable closed subsets of $X$, there is $E \in \mathcal{F}$ such that $\overline H \cap \overline E = \emptyset$. Since $X$ is $\omega$-normal, there is a real valued continuous function $h$ defined on $X$ which is equal to the constant function $n+1/2$ on $\overline E$ and is equal to $f\restriction \overline H$ on $\overline H$. It is possible now to prove that $h \in  
{\rm int}V(f,\overline H,\epsilon) \cap V(n+1/2,\overline E, 1/2)$. This means that $D_n$ is dense in $C_s(X)$. Since $C_s(X)$ is a Baire space, there is a function $g \in \bigcap _{n< \omega}D_n$. Hence, for each $n < \omega$, there exists $\overline E_n \in \mathcal{F}$ such that $g \in V(n+1/2,\overline E_n,1/2)$; that is, $n < g(x) < n+1$ for each $ x \in \overline E_n$. If we take $V_n = g^{-1}[(n,n+1)]$, then each $V_n$ is open and $E_n \subseteq V_n$. Moreover, it is easy to prove that the collection $\{V_n : n < \omega, n \,\, \text{is odd}\}$ is a discrete sequence. 
\end{proof}

\begin{corollary}
Let $X$ be a non-separable $\omega$-normal topological space such that the clousure of every countable subset is countable. If $C_s(X)$ is a Baire space, then every collection of cardinality $\omega_1$ of closed separable subsets of $X$ contains a strongly discrete sequence.
\end{corollary}

\begin{proof}
Indeed, if $\mathcal{F}$ is a family of closed separable subsets of $X$ and $|\mathcal{F}| = \omega_1$, then $\mathcal{F}$ moves off the collection of all closed separable subsets of $X$. 
\end{proof}

\begin{definition}
A topological vector space $X$ is {\it locally convex} if the origin has a local base made up of convex neighborhoods.
\end{definition}

\begin{definition}
Let $X$ be a vector space over $\mathbb{R}$ and $A \subseteq X$. We say that $A$ is {\it balanced} if, for all $\alpha \in \mathbb {R}$ with $\left \vert \alpha \right \vert \leq1 $, we have that $\alpha \cdot A \subseteq A$.
\end{definition}

\begin{definition}
A locally convex vector space $X$ is of {\it type Baire} if for any increasing sequence  $\left \{A_ {n}: n \in \mathbb {N}\right \}$ formed by convex and balanced closed sets  such that $X = \bigcup_ {n \in \mathbb{N}}A_{n}$, there exists $n \in \mathbb {N}$ such that $A_ {n}$ is a neighborhood of the origin.
\end{definition}

All Baire locally convex vector space is of type Baire.

\begin{theorem}
$C_{s}\left(X\right)$ is a locally convex vector space if and only if 
$X$ is pseudocompact.
\end{theorem}

\begin{theorem}
If $X$ is pseudocompact and $C_{s}\left(X\right)$ is of type Baire,
then $X$ is separable.
\end{theorem}

\begin{proof} 
For each $n\in\mathbb{N}$ we define 
\[
A_{n}=\left\{  f\in C_{s}\left(  X\right)  :\underset{x\in X}{\sup}\left\vert
f\left(  x\right)  \right\vert \leq n\right\}
\]
The family $\left\{A_{n}:n\in\mathbb{N}\right\}$ is an increasing sequence of convex balanced subsets of $C_s(X)$. 
Besides, each $A_{n}$ is a closed subset. Indeed, if $g\notin A_{n}$, then there exist $z\in X$ and
$s\in\mathbb{R}$ such that $\left\vert g\left(z\right) \right\vert >s>n$ and 
so we have $V\left(  g,z,s-n\right)  \subseteq C_{s}\left(
X\right)  \setminus A_{n}$. This last contention occurs since if $f\in V\left(
g,z,s-n\right)$, then%
\[
\left\vert f\left(  z\right)  \right\vert \geq\left\vert g\left(  z\right)
\right\vert -\left\vert g\left(  z\right)  -f\left(  z\right)  \right\vert
>s-\left(  s-n\right)  =n
\]
which implies that $f\notin A_{n}$.

Now, as $C_{s} \left(X \right)$ is of type Baire, there exists $n \in \mathbb{N}$ such that $0 \in int\left(A_{n}\right)$. Therefore, we can take a separable closed subspace $D \subseteq X$ and $\varepsilon> 0$ such that
\[
V\left(  0,D,\varepsilon\right)  \subseteq A_{n}.
\]

To conclude the proof, it will be shown that $D = X$. Supposing that
there is $z \in X \setminus D$, since $X$ is Tychonoff and $D$ is closed, you can
get $f\in C\left(  X\right)$ such that $f\left(  z\right)  =n+1$ and 
$f\upharpoonright_{D}\equiv0$. Hence $f\in V\left(  0,D,\varepsilon
\right)$ and therefore $\left\vert f\left(  x\right)  \right\vert \leq n$
for all $x\in X$. In particular, $f\left(  z\right)  \leq n$, which is impossible since 
$f \left(z \right) = n + 1$. This proves that $D = X$ and therefore $X$ is separable. 
\end{proof}

\begin{corollary}
If $X$ is pseudocompact and not separable, then $C_{s}\left(  X\right)$
is not of type Baire.
\end{corollary}

\section{Cellularity of $C_s(X)$ for spaces of ordinals}

 In this section we study the cellularity of $C_s([0,\gamma]\times X)$ where $[0,\gamma]$ is the collection of all ordinal numbers belonging to $\gamma$ and with the topology generated by its usual order. In particular, we prove that under certain assumptions (see Theorem \ref{knasteromegabounded}) the cellularity of this kind of spaces is always bounded by $\omega_1$. As a consequence of this we completely determine the cellularity of $C_s(X)$ where $X$ is a space of ordinals.\\
 
 The cellularity of $C_s(X)$ has already been studied by the last three authors in \cite{RTV}. In that paper, compact metrizable spaces are characterized by means of the following theorem.

 \begin{theorem}\rm{(Theorem 7.3, \cite{RTV})}\label{cellularitymetric}
 $c(C_s(X))=\omega$ if and only if $X$ is compact metrizable.     
 \end{theorem}

The previous theorem is important for  many reasons. Particularly it represents on of the many examples in which classical topological properties associated to $X$ can be studied by analyzing cardinal functions on $C_s(X)$. The next natural step is to study spaces for which $c(C_s(X))=\omega_1$. In order to do that, we will need the following lemmas.

\begin{lemma}\label{lemmaunion}Let $X$ be space and let $Y$ be a clopen subset of $X$. Then the function $H:C_s(X)\longrightarrow C_s(Y)\times C_s(X\backslash Y)$ given by $H(f)=(f|_Y,f|_{X\backslash Y})$ is an homeomorphism.
\end{lemma}

In the following proposition we make use of elementary submodels. The reader unfamiliar with this technique can learn more in \cite{Alan} or \cite{Sou}.

\begin{lemma}{\rm (Weak $\Delta$-system lemma)} Let $\mathcal{A}$ be a family of countable sets with $|\mathcal{A}|=\omega_2$. There is a countable set $R$ and a subfamily $\mathcal{B}\in[\mathcal{A}]^{\omega_2}$ with the following property:
$$\forall F,G\in \mathcal{B}\,(\,F\not=G\rightarrow F\cap G\subseteq R\,).$$ We call $\mathcal{B}$ a weak $\Delta$-system with weak root $R$.
\end{lemma}

\begin{proof}Let $\lambda$ be sufficiently large and let $M$ be an elementary submodel of $H(\lambda)$ with the following properties:
\begin{enumerate}
\item$|\mathcal{A}|=\omega_1,$

\item$\mathcal{A}\in M$,

\item For every $X\in [M]^{\omega}$ there is $R\in M\cap [M]^{\omega}$ for which $X\subseteq R.$ In particular, this implies that $X\subseteq M$ whenever $X\in M$ and $|X|\leq \omega_1.$
\end{enumerate}

Since the cardinality of $\mathcal{A}$ is greater than $\omega_1$, there is $A\in \mathcal{A}\backslash M.$ Note that $A\cap M$ is countable, so there is $R\in M\cap [M]^{\omega}$ for which $A\cap M\subseteq R$.  Zorn's lemma implies that there is a subfamily of $\mathcal{A}$ maximal with respect to the property of being a Weak $\Delta$-system with weak root $R$. Using elementarity, we can take $\mathcal{B}\in M$ as a witness of the previous fact. To finish, it is enough to prove that $|\mathcal{B}|=\omega_2.$ Suppose towards a contradiction that this is not the case. By property $(3)$, it must happen that $\mathcal{B}\subseteq M$. Take an arbitrary $B\in \mathcal{B}$, since $B$ is countable (remember that $B\in \mathcal{A}$) it is also true that $B\subseteq M$. Thus, $B\cap A=B\cap M\cap A\subseteq B\cap R\subseteq R$. In this way, $\mathcal{B}\cup\{A\}$ is a weak $\Delta$-system with weak root $R$. Consequently, $\mathcal{B}\cup\{A\}=\mathcal{B}$ which means that $A\in M$. This is a contradiction so we are done.
\end{proof}

\begin{definition}Let $X$ be a space and $\kappa$ be a cardinal number. We say that $X$ is $\kappa$-Knaster
 if each family of $\mathcal{U}$ of non-empty open sets with $|\mathcal{U}|=\kappa$ there is $\mathcal{V}\in [\mathcal{U}]^{\kappa}$ such that $U\cap V\not=\emptyset$ for any $U,V\in\mathcal{V}.$
\end{definition}

The following proposition resumes all the properties which related to $\kappa$-Knaster spaces that we will use.

\begin{proposition}\label{knasterproperties}Let $X$ and $Y$ be two spaces and let $\kappa$ be a cardinal number. Then:

\begin{enumerate}
\item If $d(X)\leq \kappa$, then $X$ is $\kappa$-Knaster.

\item If $X$ and $Y$ are both $\kappa$-Knaster, then $X\times Y$ is $\kappa$-Knaster.

\item If $X$ is $\kappa$-Knaster and $Y$ is a continuous image of $X$, then $Y$ is also $\kappa$-Knaster.

\item If $X$ is $\kappa$-Knaster, then $X$ has no cellular family of size $\kappa$. In particular, if $\kappa$ is a successor cardinal then $c(X)<\kappa.$
\end{enumerate}
\end{proposition}

As a corollary of Lemma \ref{lemmaunion} and the point $(2)$ in Proposition \ref{knasterproperties} we have the following.

\begin{corollary}\label{knasterunion} Let $\kappa$ be a cardinal number, $X$ be a space and $Y\subseteq X$ be a clopen subspace of $X$. If $C_s(Y)$ and $C_s(X\backslash Y)$ are both $\kappa$-Knaster then $C_s(X)$ is $\kappa$-Knaster.
\end{corollary}

\begin{theorem}\label{knasteromegabounded} Let $X$ be an $\omega$-bounded space for which $C_s(X)$ is $\omega_2$-Knaster. Then $C_s([0,\gamma]\times X))$ is $\omega_2$-Knaster for any non-empty ordinal $\gamma$. In particular, $c(C_s([0,\gamma]\times X))\leq \omega_1.$
\end{theorem}

\begin{proof}Fix $X$ as in the hypotheses. We will prove the theorem by induction over $\gamma$. If $\gamma=1$ the result is immediate. For the inductive step of the proof, we first 
deal with the case where $\gamma=\gamma'+1$ and the theorem holds for $\gamma'$.  it suffices to note that $[0,\gamma]\times X$ can be partitioned into the clopen sets $[0,\gamma']\times X$ and $\{\gamma\}\times X$. Thus, the conclusion holds due to the induction hypothesis and Corollary \ref{knasterunion}. Finally, we deal with the case where $\gamma$ is a limit ordinal and the result hold for each $\alpha\in[0,\gamma]$. We start by taking a family of $\omega_2$ open sets in $C_s([0,\gamma]\times X)$. Without loss of generality that family can be written as $\langle V(f_\xi,F_\xi,\epsilon)\rangle_{\xi\in \omega_2}$ where the following conditions hold for each $\xi\in\omega_2$:

\begin{itemize}
\item $f_\xi\in C_s([0,\gamma]\times X)$,

\item $F_\xi=G_\xi\times A_\xi$ where $G_\xi$ and $A_\xi$ are closed separable subsets of $[0,\gamma]$ and $X$ respectively and $\gamma\in G_\xi$. In particular,  both $A_\xi$ and $G_\xi$ are compact $G_\xi$ is countable.
\end{itemize}

Our goal is to find distinct $\xi,\xi'\in \omega_2$ for which $V(f_{\xi},F_\xi,\epsilon_\xi)\cap V(f_{\xi'},F_{\xi'},\epsilon_\xi')\not=\emptyset.$\\ Given $\xi\in \omega_2$, let $h_\xi:X\longrightarrow \mathbb{R}$ be defined as:
$$h_\xi(x)=f_\xi(\gamma,x).$$ Since $C_s(X)$ is $\omega_2$-Knaster, we can find $\mathcal{A}\in[\omega_2]^{\omega_2}$ with the property that  any two elements of the family $\langle V(h_\xi,A_\xi, \frac{\epsilon}{2})\rangle_{\xi\in \mathcal{A}}$ have non-empty intersection. Momentarily let us fix $\xi \in \mathcal{A}$. As $f_\xi$ is continuous, for each $x\in A_\xi$ there are $\beta_x<\gamma$ and an open neighborhood $U_x\subseteq X$ of $x$  for which $f_\xi\big[ [\beta_x,\gamma]\times U_x\big]\subseteq (x-\frac{\epsilon}{4},x+\frac{\epsilon}{4})$. Note that $|f_\xi(\gamma,y)-f_\xi(\alpha,y)|<\frac{\epsilon}{2}$ for each $y\in A$ and $\alpha\in [\beta_x,\gamma]$. Now,  $A_\xi$ is compact so there is a finite $S
\subseteq A_\xi$ for which $\{ U_x\,|\,x\in S\}$ covers $A_\xi$. Consider $\beta_\xi=\max\{\beta_x\,|\,x\in S\}$. By a previous observation, $$|f_\xi(\gamma,x)-f_\xi(\alpha,x)|<\frac{\epsilon}{2}$$ for each $x\in A_\xi$ and $\alpha\in [\beta_\xi,\gamma]$.\\

The rest of the proof is handled in two cases.\\

\noindent {\bf First case.} If $cof(\gamma)=\omega$. In this case we can take a countable set $M\subseteq \gamma$ with $\sup M=\gamma$. For every $\xi\in \mathcal{A}$ there is $\beta \in M$ for which $\beta_\xi<\beta$. Since $\omega_2$ is regular and $|\mathcal{A}|=\omega_2$ there must be $\mathcal{B}\in [\mathcal{A}]^{\omega_2}$ and $\beta\in M$ such that $\beta_\xi<\beta$ for any $\xi\in\mathcal{B}.$ Now consider the collection $\langle V(f_\xi|_{[0,\beta]},F_\xi\cap \big([0,\beta]\times X\big),\epsilon)\rangle_{\xi\in \mathcal{B}}$ of open sets in $C_s([0,\beta]\times X)$. By induction hypothesis there is $\mathcal{C}\in [\mathcal{B}]^{\omega_2}$ for which any two elements of $\langle V(f_\xi|_{[0,\beta]\times X},F_\xi\cap \big([0,\beta]\times X\big),\epsilon)\rangle_{\xi\in \mathcal{C}}$ have non-empty intersection. We claim that the same is true for $\langle V(f_\xi, F_\xi, \epsilon)\rangle_{\xi\in \mathcal{C}}$. For this purpose, let $\xi$ and $\xi'$ be arbitrary elements of $\mathcal{C}$.  Let $g\in C_s([0,\beta]\times X)$ and $h\in C_s(X)$ be such that:

\begin{itemize}
\item $g\in V(f_\xi|_{[0,\beta]\times X},F_\xi\cap\big([0,\beta]\times X\big),\epsilon)\cap V(f_{\xi'}|_{[0,\beta]\times X},F_{\xi'}\cap\big([0,\beta]\times X\big),\epsilon)$,

\item $h\in V(h_\xi,A_\xi,\frac{\epsilon}{2})\cap V(h_{\xi'},A_{\xi'}\frac{\epsilon}{2}).$
\end{itemize}

As $[0,\beta]\times X$ and $[\beta+1,\gamma]\times X$ form a partition of $[0,\gamma]\times X$ into clopen sets, then the function $f:[0,\gamma]\times X\longrightarrow \mathbb{R}$ given as:$$f(\alpha,x)=\begin{cases} g(\alpha,x)&\textit{if } \alpha\leq \beta\\
h(x)&\textit{ if }\alpha>\beta
\end{cases}$$ belongs to $ C_s([0,\gamma]\times X).$ To finish, let $(\alpha,x)\in F_\xi$. If $\alpha\leq \beta$ then $|f(\alpha,x)-f_\xi(\alpha,x)|=|g(\alpha,x)-f_\xi(\alpha, x)|<\epsilon$. On the other side, if $\alpha>\beta$ then \begin{align*} 
|f(\alpha,x)-f_\xi(\alpha,x)|=&|h(x)-f_\xi(\alpha,x)|\\= & |h(x)-f_\xi(\gamma,x)+f_\xi(\gamma,x)-f_\xi(\alpha,x)|\\ \leq & |h(x)-h_\xi(x)|+|f_\xi(\gamma,x)-f_\xi(\alpha,x)|\\ < &\frac{\epsilon}{2}+\frac{\epsilon}{2}=\epsilon.
\end{align*}
In this way, we showed that $\hat{h}\in V(f_\xi, F_\xi, \epsilon)$. Analogously, we can show that $\hat{h}\in V(f_{\xi'},F_{\xi'},\epsilon)$ so we are done.\\

\noindent {\bf Second case.} If $cof(\gamma)>\omega$. By the weak $\Delta$-system lemma, there are $\mathcal{B}\in [\mathcal{A}]^{\omega_2}$ and a countable $R\subseteq [0,\gamma]$ for which $\{G_\alpha\}_{\alpha\in \mathcal{B}}$ is a weak $\Delta$-system with weak root $R$. Since the closure of countable set os $[0,\gamma]$ is again countable, we can suppose without loss of generality that $R$ is closed and has at least two points. A $\gamma$ has uncountable cofinality and $R$ is countable, it must happen that $\gamma$ is an isolated point of $R$. In this way, $\beta:=\max(R\backslash \{\gamma\})$ is well defined and is an ordinal strictly smaller than $\gamma$. As in the previouse case, we use the induction hypothesis to find $\mathcal{C}\in [\mathcal{B}]^{\omega_2}$ so that any two elements of $\langle V(f_\xi|_{[0,\beta]\times X},F_\xi\cap \big([0,\beta]\times X\big),\epsilon)\rangle_{\xi\in \mathcal{C}}$ have non-empty 
intersection. As in the previous case, we claim that the same holds for $\langle V(f_\xi, F_\xi,\epsilon)\rangle_{\xi \mathcal{C}}.$ Indeed, let $\xi$ and $\xi'$ be arbitrary elements of $\mathcal{C}$. Note that $\gamma$ is an isolated point of both $G_\xi$ and $G_{\xi'}$. In this way, $H_\xi:=\big( G_\xi\cap [\beta+1,\gamma)\big)\times A_\xi$ and $H_{\xi'}\big(G_{\xi'}\cap [\beta+1,\gamma)\big)\times A_{\xi'} $ are both closed subsets of $[0,\gamma]\times X$. Furthermore, \begin{align*}H_{\xi} \cap H_{\xi'}\subseteq &\big( G_\xi\cap G_{\xi'}\cap[\beta,\gamma)\big)\times X\\ \subseteq &\big( R\cap [\beta,\gamma)\big)\times X=\emptyset.
\end{align*}
Let $g\in C_s([0,\beta]\times X)$ and $h\in C_s(X)$ be such that:\begin{itemize}
\item $g\in V(f_\xi|_{[0,\beta]\times X},F_\xi\cap\big([0,\beta]\times X\big),\epsilon)\cap V(f_{\xi'}|_{[0,\beta]\times X},F_{\xi'}\cap\big([0,\beta]\times X\big),\epsilon)$,
\item $h\in V(h_\xi,A_\xi,\frac{\epsilon}{2})\cap V(h_{\xi'},A_{\xi'}\frac{\epsilon}{2}).$
\end{itemize}
As $H_\xi, H_{\xi'}$, $\big(\{\gamma\}\times X\big)\cap \big(F_\xi\cup F_{\xi'}\big)$ and $\big([0,\beta]\times X\big)\cap \big( F_\xi\cup F_{\xi'})$ are pairwise disjoint closed sets, the function $\hat{h}:F_\xi\cup F_{\xi'}\longrightarrow \mathbb{R}$ defined as:
$$\hat{h}(\alpha,x)=\begin{cases}g(\alpha,x)&\textit{if }\alpha\leq \beta\\ f_\xi(\alpha,x)&\textit{if }(\alpha,x)\in H_\xi\\
f_{\xi'}(\alpha,x)&\textit{if }(\alpha,x)\in H_{\xi'}\\
h(x)&\textit{if }\alpha=\gamma
\end{cases}$$
Is continuous. Furthermore , the domain of $\hat{h}$ is compact which means that we can use the Tietze-Urysohn extension Theorem to find $\hat{\hat{h}}\in C_s([0,\gamma]\times X) $ which extends $\hat{h}$. By arguing in a similar manner as in the case where $cof(\gamma)=\omega$, we can conclude that $\hat{\hat{h}}\in V(f_\xi,F_\xi,\epsilon)\cap V(f_{\xi'},F_{\xi'},\epsilon)$. This finishes the proof.
\end{proof}

In Theorem 5.6 of \cite{PichTamVill} it is proved that $d(C_s(X))=\omega$ whenever $X$ is a compact metrizable space. Hence, by the point $(1)$ of Proposition \ref{knasterproperties} we have the following corollary. 

\begin{corollary}\label{coro:prod} Let $X$ be a compact metrizable space. Then $C_s([0,\gamma]\times X)$ is $\omega_2$-Knaster for each ordinal $\gamma.$ In particular, $c(C_s([0,\gamma] \times X)) = \omega$ if $\gamma$ is countable, and $c(C_s([0,\gamma] \times X)) = \omega_1$ if $\gamma$ is uncountable. 
\end{corollary}

\begin{proof}
The first assertion is a consequence of Theorem \ref{knasteromegabounded}. The second assertion follows from Theorem \ref{cellularitymetric}.
\end{proof}

\begin{corollary}\label{coro:prod}
Let $n\in \omega$ and $\gamma_0,\dots,\gamma_n$ be an arbitrary collection of ordinals. Then $C_s(\prod_{i\leq n}[0,\gamma_i])$ is $\omega_2$-Knaster. In particular, $c(C_s(\prod_{i\leq n}[0,\gamma_i])) = \omega$ if for each $i \leq n$, $\gamma_i$ is countable, and $c(C_s(\prod_{i\leq n}[0,\gamma_i])) = \omega_1$ if there is $i \leq n$ such that $\gamma_i$ is uncountable. 
\end{corollary}

\begin{proof} The proof is carried by induction over $n$. For the case where $n=0$, let $X$ be a space with only one point. $X$ trivially satisfies the conditions of Theorem \ref{knasteromegabounded}. This means that $C_s([0,\gamma_0]\times X)$ is $\omega_2$-Knaster. As $C_s([0,\gamma_0]\times X)$ is homeomorphic to $C_s([0,\gamma_0])$, we are done. Now suppose that we have proved the corollary for some $n$ and let $\gamma_0\dots,\gamma_{n+1}$ be distinct ordinals. By induction hypothesis we have that $C_s(\prod_{1\leq i\leq n}[0,\gamma_i])$ is $\omega_2$-Knaster. It is well known that $\prod_{1\leq i\leq n}[0,\gamma_i]$ is $\omega$-bounded, so by Theorem \ref{knasteromegabounded} we have that  $C_s([0,\gamma_0]\times \prod_{1\leq i\leq n+1}[0,\gamma_i])=C_s(\prod_{i\leq n+1}[0,\gamma_i])$ is $\omega_2$-Knaster. This finishes the proof of the first part. The second part can be proved in a symilar way than the proof of Corollary \ref{coro:prod}.
\end{proof}

\begin{corollary}
Let $\gamma$ be an ordinal, then $c(C_s([0,\gamma])) = \omega$ if $\gamma$ is countable, and $c(C_s([0,\gamma])) = \omega_1$ if $\gamma$ is uncountable.
\end{corollary}

\begin{proposition}\label{pseudocompactcontinuous} Let $\gamma$ be an ordinal and let $X$ be a countably compact subspace of $[0,\gamma]$. Then $C_s(X)$ is a continuous image of $C_s(cl_{[0,\gamma]}{X}).$
\end{proposition}

\begin{proof} Let $p:C_s(cl_{[0,\gamma]}{X})\longrightarrow C_s(X)$ given by $p(f)=f|_X$. Clearly $p$ is continuous (and injective). We will show that it is surjective. For this, let $f\in C_s(X)$ and take an arbitrary $\beta\in cl_{[0,\gamma]}{X}\backslash X$. Since $X$ is pseudocompact, it must be the case that $cof(\beta)>\omega$. If this was not the case, there would be a sequence in $X$ converging to $\beta$. Hence, that sequence would form a countable closed discrete set in $X$ which is impossible. As $cof(\beta)>\omega$ we know there is $\alpha<\beta$ for which $f|_{X\cap[\alpha,\beta)}$ is constant with value $r_\beta$. It is straightforward that the function $\hat{f}:cl_{[0,\gamma)}{X}\longrightarrow \mathbb{R}$ given by:$$\hat{f}(\beta)=\begin{cases}f(\beta)&\textit{if }\beta\in X\\
r_\beta &\textit{otherwise}
    
\end{cases}$$
belongs to $C_s(cl_{[0,\gamma)}{X})$ and $p(\overline{f})=f$.
\end{proof}

As a direct consequence of the previous proposition and by the point $(3)$ of Proposition \ref{knasterproperties} we have that $c(C_s(X))\leq\omega_1$ whenever $X$ is a countably compact space of ordinals. So, we have:

\begin{corollary}\label{coro:subspaces}
If $X$ is a countably compact subspace of $[0,\gamma]$, then 
\begin{enumerate}
\item $c(C_s(X)) = \omega$ if and only if $X$ is compact and countable.

\item $c(C_s(X)) = \omega_1$ if and only if $X$ is uncountable.
\end{enumerate}

In particular, if $cf(\gamma) > \omega$, then  
$c(C_s([0,\gamma)) = \omega_1$.
\end{corollary}

In Corollary 5.11 of \cite{PichTamVill} it is proved that $c(C_s(X))\geq 2^\omega$ whenever $X$ is not pseudocompact. Since pseudocompactness and countably compactness coincide for spaces of ordinals, we have the following corollary.

\begin{corollary}\label{coro:CH}
Under the negation of the continuum hypothesis, if $X$ is a non-countable space of ordinals, then $X$ is pseudocompact if and only if $c(C_s(X))=\omega_1$. 
\end{corollary}

It is still true that if $X$ is a countably compact space of finite products of ordinals then $c(C_s(X)) \leq \omega_1$. To prove this, we will need a different aproach from the one used in Corollary \ref{coro:subspaces}.

\begin{proposition} 
Let $\gamma$ be an ordinal and $n \in \omega$. If $X \subseteq [0,\gamma]^n$
is countably compact, then $C_s(X)$ is $\omega_2$-Knaster.
\end{proposition}

\begin{proof} 
Let $\langle V_\alpha \rangle_{\alpha \in \omega_2}$ be a family of nonempty open sets in $C_s(X)$. Without loss of generality,
there is $\epsilon_\alpha > 0$ such that each $V_\alpha$ is of the form $V(f_\alpha; F_\alpha;\epsilon_\alpha)$ for some $f_\alpha \in C_s(X)$ and $F_\alpha$ closed
and separable in $X$. Let us fix $\alpha \in \omega_2$. Since $F_\alpha$ is closed in $X$, then $F_\alpha$ is countably compact. Furthermore, since $F_\alpha$ is separable and it is a subspace of $[0,\gamma]^n$ it is nessesarily countable. This
means that $F_\alpha$ is compact, which in particular means that it is closed in 
$[0,\gamma]^n$. Hence, by the Tietze-Urysohn extension Theorem there is $g_\alpha\in C_s([0,\gamma]^n)$ extending $f_\alpha \restriction F_\alpha$. As 
$C_s([0,\gamma]^n)$ is $\omega_2$-Knaster, there is $\mathcal{A} \in [\omega_2]^{\omega_2}$ such that $V(g_\alpha;F_\alpha;\epsilon_\alpha) \cap V(g_\beta;F_\beta;\epsilon_\beta) \not= \emptyset$ for all $\alpha, \beta \in \mathcal{A}$. To finish the proof, just note that if $h \in V(g_\alpha;F_\alpha;\epsilon_\alpha) \cap V(g_\beta;F_\beta;\epsilon_\beta)$, then $h \restriction X \in V_\alpha \cap V_\beta$.
\end{proof}

Unfortunately there are pseudocompact subpaces of finite products of ordinals which fail to be countably compact. This means that we cannot use the same argument as before to generalize Corollary \ref{coro:CH}. This leaves us with the following problem.

\begin{problem}
Is it consistent with the negation of the continuum hypothesis that there is a pseudocompact subspace of finite product of ordinals, say X, such that $c(C_s(X)) \geq 2^\omega$?
\end{problem}

As a corollary of Theorem 7.18 of \cite{RTV} we have that $c(C_s([0;\kappa]^\kappa)) = \kappa$ for each cardinal $\kappa$. On the
other side $c(C_s[0,\kappa]^n) = \omega_1$ for each $n \in \omega$, due to Corollary \ref{coro:prod}. Thus, it is natural to ask about the
situation in intermediate products. In particular, we ask the following.

\begin{problem}
Which are the values of $c(C_s([0,\omega_2]^\omega))$ and $c(C_s([0,\omega_2]^{\omega_1}))$?
\end{problem}

The reader interested in knowing more about topological properties over products of ordinals may wish to search for \cite{KS} and \cite{KOT}.

\end{document}